\documentclass[12pt]{article}
\usepackage{latexsym,amsfonts,amssymb,mathrsfs,float}
\usepackage{dsfont}
\usepackage{amsthm}
\usepackage[numbers]{natbib}
\usepackage{amsmath}
\usepackage{indentfirst}
\usepackage{enumerate}

\usepackage[dvips]{color}
\usepackage{colordvi,multicol}
\allowdisplaybreaks

\newtheorem{thm}{Theorem}[section]

\newtheorem{lem}[thm]{Lemma}

\newtheorem{rem}[thm]{Remark}

\numberwithin{equation}{section}
\begin{document}

\title{\bf A study on the negative binomial distribution motivated by Chv\'{a}tal's theorem}

 \author{ Zheng-Yan Guo, Ze-Yu Tao, Ze-Chun Hu\thanks{Corresponding
author: zchu@scu.edu.cn}\\ \\
 {\small College of Mathematics, Sichuan  University,
 Chengdu 610065, China}}

 \date{}
\maketitle

\begin{abstract}

  Let $B(n,p)$ denote a binomial random variable with parameters $n$ and $p$.  Chv\'{a}tal's theorem  says that for any fixed $n\geq 2$, as $m$ ranges over $\{0,\ldots,n\}$, the probability $q_m:=P(B(n,m/n)\leq m)$ is the smallest when  $m$ is closest to $\frac{2n}{3}$.  Motivated by this theorem, in this note  we consider the infimum value of the probability $P(X\leq  E[X])$, where   $X$ is  a negative binomial random variable. As a consequence, we give an affirmative answer to the conjecture posed in [Statistics and Probability Letters, 200 (2023) 109871].

\end{abstract}



\section{Introduction and main result}\setcounter{equation}{0}

Let $B(n,p)$ denote a binomial random variable with parameters $n$ and $p$. Chv\'{a}tal's theorem  says that for any fixed $n\geq 2$, as $m$ ranges over $\{0,\ldots,n\}$, the probability $q_m:=P(B(n,m/n)\leq m)$ is the smallest when  $m$ is closest to $\frac{2n}{3}$.

Chv\'{a}tal's theorem has  applications in machine learning (see \cite{Do18}, \cite{GM14} and the references therein). As to the history of Chv\'{a}tal's theorem, we refer to   Barabesi et al. \cite{BPR21}, Janson \cite{Ja21} and Sun \cite{Su21}.

Motivated by Chv\'{a}tal's theorem, the infimum value problem of the probability
$P(X\leq E[X])$ has been studied in a series paper recently, where $X$ is a random variable and $E[X]$ is its expectation. In the following, we list the distributions which have been considered in this aspect.

\begin{itemize}

\item The Poisson distribution, the geometric distribution and the Pascal distribution (Li et al. \cite{LXH23});

\item The Gamma distribution (Sun et al. \cite{SHS23});

\item  The Weibull distribution and the Pareto distribution (Li et al. \cite{LHZ23});

\item The inverse Gaussian, log-normal, Gumbel and logistic distributions (Hu et al. \cite{HLZZ23}).
\end{itemize}

As to the Pascal distribution, \cite{LXH23} gave only a partial result. Let $B^*(r,p)$ denote a Pascal random variable with parameters $r$ ($r$ is a positive integer) and $p\, (0<p\leq1)$ such that
\begin{eqnarray*}
P(B^*(r,p)=j)=
\left(
\begin{array}{l}
j-1\\
r-1
\end{array}
\right)(1-p)^{j-r}p^r,\ j=r,r+1,\ldots.
\end{eqnarray*}
We know that the expectation  of $B^*(r,p)$ is $r/p$.  An anonymous referee to the first version of \cite{LXH23} posed the following {\bf conjecture}:
\begin{eqnarray}\label{1.1}
\inf_{0<p\leq 1}P\left(B^*(r,p)\leq \frac{r}{p}\right)=\left(\frac{r}{r+1}\right)^r.
\end{eqnarray}
\cite{LXH23} proved that (\ref{1.1}) holds for $r=1,2,3$.

The motivation of this note is to give an affirmative answer to (\ref{1.1}) for any positive integer $r$. In fact, we will build the corresponding result for the negative binomial distribution, which extends the Pascal distribution in some sense.

Let $NB(r,p)$ denote a negative binomial random variable with parameters $r$ ($r$ is a positive real number) and $p\, (0<p\leq1)$ such that
 \begin{eqnarray}\label{1.1-b}
P(NB(r,p)=l)=\binom{-r}{l}p^r(-q)^l=\binom{r+l-1}{l}p^rq^{l},\ l=0,1,2,\ldots,
\end{eqnarray}
where $q=1-p$. We know that the expectation of $NB(r,p)$ is $\frac{rq}{p}$.

The main result of this note is

\begin{thm}\label{thm}
For any $r>0$, it holds that
\begin{eqnarray}\label{1.2}
\inf_{0<p\leq 1}P\left(NB(r,p)\leq \frac{r(1-p)}{p}\right)=\left(\frac{r}{r+1}\right)^r.
\end{eqnarray}
\end{thm}

\begin{rem}
Note that if $r$ is a positive integer, then $NB(r,p)+r$ is a Pascal random variable with parameters $r$  and $p$, and thus
$$
P\left(NB(r,p)\leq \frac{r(1-p)}{p}\right)=P\left(NB(r,p)+r\leq \frac{r}{p}\right)=P\left(B^*(r,p)\leq \frac{r}{p}\right).
$$
It follows that Theorem \ref{thm} gives an affirmative answer to Conjecture (\ref{1.1}) for any positive integer $r$.
\end{rem}

\section{Proof of Theorem \ref{thm}}\setcounter{equation}{0}

In this section, we will give the proof of Theorem \ref{thm}. At first, we prove two lemmas.  For convenience, we recall the definitions of the Beta function $\mathbf{B}(a,b)$ and the Gamma function $\Gamma(s)$ as follows:
\begin{align*}
&\mathbf{B}(a,b):=\int_0^1x^{a-1}(1-x)^{b-1}dx,\ a>0, b>0;\\
&\Gamma(s):=\int_0^{\infty}x^{s-1}e^{-x}dx,\ s>0.
\end{align*}
It is well known that
\begin{align}
&\mathbf{B}(a,b)=\frac{\Gamma(a)\Gamma(b)}{\Gamma(a+b)},\ \forall a>0,b>0;\label{2.0}\\
&\Gamma(s+1)=s\Gamma(s),\ \forall s>0.\label{2.0-b}
\end{align}
By (\ref{2.0-b}), for any positive real number $s$,  we can define
\begin{align}\label{2.0-c}
s!:=\Gamma(s+1).
\end{align}
If $s=n$ for a nonnegative integer $n$, then the above definition  is just the usual definition, i.e.
\begin{eqnarray*}
n!=\left\{
\begin{array}{cl}
1,&\mbox{if}\ n=0;\\
\Pi_{k=1}^nk,& \mbox{if}\  n=1,2,\ldots.
\end{array}
\right.
\end{eqnarray*}

\smallskip

\begin{lem}\label{lem-2.1} For any $n\in \mathbb{N}$ and any $r\in (0,\infty)$, it holds that
 \begin{align}\label{2.1}
\left(\frac{r}{n+r+1}\right)^r\sum_{k=0}^{n}\binom{k+r-1}{k}\left(1-\frac{r}{n+r+1}\right)^k
&=\frac{1}{\mathbf{B}(n+1,r)} \int_0^{\frac{r}{n+1+r}} x^{r-1}(1-x)^ndx.
\end{align}
\end{lem}

\noindent {\bf Proof.} Let $\lambda:=\frac{r}{n+r+1}$. Then (\ref{2.1}) becomes
\begin{align}\label{2.2}
\lambda^r\sum_{k=0}^{n}\binom{k+r-1}{k}\left(1-\lambda\right)^k
=\frac{1}{\mathbf{B}(n+1,r)}\int_0^{\lambda}x^{r-1}(1-x)^ndx.
\end{align}

By the binomial theorem, we have
$$
(1-x)^n=(-x)^n+\binom{n}{n-1}(-x)^{n-1}+\cdots+\binom{n}{1}(-x) +1.
$$
It follows that
\begin{align*}
x^{r-1}(1-x)^n=(-1)^nx^{n+r-1}+\binom{n}{n-1}(-1)^{n-1}x^{n+r-2}+\cdots+\binom{n}{1}
(-1)x^{r} +x^{r-1},
\end{align*}
and thus
\begin{align}\label{2.4}
\int_0^{\lambda}x^{r-1}(1-x)^ndx=
\frac{(-1)^n}{n+r}\lambda^{n+r}+\frac{(-1)^{n-1}}{n+r-1}\binom{n}{n-1}\lambda^{n+r-1}+
\cdots+\frac{(-1)}{r+1}\binom{n}{1}\lambda^{r+1}+\frac{1}{r}\lambda^r.
\end{align}

By  (\ref{2.4}), we know that (\ref{2.2}) becomes
\begin{align}\label{2.5}
&\lambda^r\sum_{k=0}^{n}\binom{k+r-1}{k}\left(1-\lambda\right)^k\nonumber\\
&=
\frac{1}{\mathbf{B}(n+1,r)}\left[\frac{(-1)^n}{n+r}\lambda^{n+r}+\frac{(-1)^{n-1}}{n+r-1}\binom{n}{n-1}\lambda^{n+r-1}+
\cdots+\frac{(-1)}{r+1}\binom{n}{1}\lambda^{r+1}+\frac{1}{r}\lambda^r\right].
\end{align}
Obviously, the left side of the above equality can also be regarded as a polynomial of $\lambda$. Hence, in order to prove (\ref{2.5}),  we  need only to check the coefficients for each degree of $\lambda$. We decompose the check into two steps.

{\it Step 1.} To prove that the coefficients of $\lambda^r$ on both sides of (\ref{2.5}) are equal.

By (\ref{2.5}), (\ref{2.0-b}) and (\ref{2.0-c}), we know that  the coefficient of $\lambda^r$ on the right side of (\ref{2.5}) is
\begin{align*}
\frac{1}{\mathbf{B}(n+1,r)r}=\frac{\Gamma(n+r+1)}{\Gamma(n+1)\Gamma(r)r}=\frac{(n+r)!}{n!(r-1)!r}=\binom{n+r}{r}.
\end{align*}
The coefficient of $\lambda^r$ on the left side of (\ref{2.5}) is
\begin{align*}
\binom{n+r-1}{n}+\binom{n+r-2}{n-1}+\cdots+\binom{r-1}{0}.
\end{align*}
By the combinatorial identity that $\binom{m}{k}=\binom{m-1}{k}+\binom{m-1}{k-1}$, we have \begin{align*}
\binom{n+r}{n}&=\binom{n+r-1}{n}+\binom{n+r-1}{n-1}\\
&=\binom{n+r-1}{n}+\binom{n+r-2}{n-1}+\binom{n+r-2}{n-2}\\
&=\cdots\\
&=\binom{n+r-1}{n}+\binom{n+r-2}{n-1}+\cdots+\binom{r}{1}+\binom{r}{0}\\
&=\binom{n+r-1}{n}+\binom{n+r-2}{n-1}+\cdots+\binom{r}{1}+\binom{r-1}{0}.
\end{align*}
Hence the coefficients of $\lambda^r$ on both sides of (\ref{2.5}) are equal.

{\it Step 2.} To prove that the coefficients of $\lambda^{r+m}$ on both sides of (\ref{2.5}) are equal  for any $1\leq m \leq n$.

The coefficient on the left side is $(-1)^m\sum_{i=m}^{n}\binom{i+r-1}{i}\binom{i}{m}$, and the one on the right side is $=(-1)^m\frac{1}{r+m}\binom{n}{m}\frac{1}{\mathbf{B}(n+1,r)}$. So it is enough to prove
\begin{align}\label{2.6}
\sum_{i=m}^{n}\binom{i+r-1}{i}\binom{i}{m}=\frac{1}{r+m}\binom{n}{m}\frac{1}{\mathbf{B}(n+1,r)}.
\end{align}

For $i=m,\ldots,n$, define
\begin{eqnarray}
&&a_i:=\binom{i+r-1}{i}\binom{i}{m},\label{2.6-a}\\
&&b_i:=\left\{
\begin{array}{cl}
\frac{1}{r+m}\binom{m}{m}\frac{1}{\mathbf{B}(m+1,r)},& \mbox{if}\ i=m,\\
\frac{1}{r+m}\left[\binom{i}{m}\frac{1}{\mathbf{B}(i+1,r)}-
\binom{i-1}{m}\frac{1}{\mathbf{B}(i,r)}\right],& \mbox{if}\ i=m+1,\ldots,n.\label{2.6-b}
\end{array}
\right.
\end{eqnarray}
Then (\ref{2.6}) becomes
$$
\sum_{i=m}^na_i=\sum_{i=m}^nb_i.
$$
Thus it is enough to show that for any $i=m,\ldots,n$, we have $a_i=b_i$.

For $i=m$, we have
\begin{align*}
b_m&=\frac{1}{r+m}\binom{m}{m}\frac{1}{\mathbf{B}(m+1,r)}=\frac{1}{r+m}\cdot\frac{(m+r)!}{m!(r-1)!}\\
&=\frac{(m+r-1)!}{m!(r-1)!}=\binom{m+r-1}{m}=a_m.
\end{align*}

For $i=m+1,\ldots,n$, by (\ref{2.6-b}), (\ref{2.0}), (\ref{2.0-c})  and (\ref{2.6-a}),  we have
\begin{align*}
b_i&=\frac{1}{r+m}\left[\binom{i}{m}\frac{1}{\mathbf{B}(i+1,r)}-
\binom{i-1}{m}\frac{1}{\mathbf{B}(i,r)}\right]\\
&=\frac{1}{r+m}\left[\frac{i!}{m!(i-m)!}\cdot\frac{(i+r)!}{i!(r-1)!}-
\frac{(i-1)!}{m!(i-m-1)!}\cdot\frac{(i+r-1)!}{(i-1)!(r-1)!}\right]\\
&=\frac{1}{r+m}\cdot\frac{1}{m!(r-1)!}\left[\frac{(i+r)!}{(i-m)!}-
\frac{(i+r-1)!}{(i-m-1)!}\right]\\
&=\frac{1}{r+m}\cdot\frac{1}{m!(r-1)!}\cdot \frac{(i+r-1)![(i+r)-(i-m)]}{(i-m)!}\\
&=\frac{(i+r-1)!}{(r-1)!i!}\cdot \frac{i!}{m!(i-m)!}\\
&=\binom{i+r-1}{i}\binom{i}{m}\\
&=a_i.
\end{align*}

Hence (\ref{2.5}) holds and the proof is complete.\hfill\fbox

\bigskip


\begin{lem}\label{lem-2.2}
For any $n\in \mathbb{N}$ and any $r\in (0,\infty)$, it holds that
\begin{align}\label{2.8}
\int_{\frac{n+1+r}{r}}^{+\infty} (u-1)^n u^{-n-r-2}\left(\frac{ru}{n+1+r}-1\right)du=\frac{r^r(n+1)^{n+1}}{(n+r+1)^{n+r+2}}.
\end{align}
\end{lem}

\noindent {\bf Proof.} Let $\mu:=\frac{r}{n+r+1}$. Then we have
\begin{align}\label{2.9}
\frac{r^r(n+1)^{n+1}}{(n+r+1)^{n+r+2}}=\frac{1}{r} \mu^{r+1} (1-\mu)^{n+1}.
\end{align}
By the change of variable $x=\frac{1}{u}$, we get
\begin{align}\label{2.10}
&\int_{\frac{n+1+r}{r}}^{+\infty} (u-1)^n u^{-n-r-2}\left(\frac{ru}{n+1+r}-1\right)du\nonumber\\
&=\int_{\frac{1}{\mu}}^{+\infty} (u-1)^n u^{-n-r-2}\left(\mu u-1\right)du\nonumber\\
&=\int_0^{\mu} (1-x)^nx^r\left(\frac{\mu}{x}-1\right)dx\nonumber\\
&=\mu\int_0^{\mu} (1-x)^nx^{r-1}dx-\int_0^{\mu} (1-x)^nx^rdx\nonumber\\
&=\mu\int_0^{\mu} \left[\sum_{k=0}^n(-1)^k\binom{n}{k}x^k\right]x^{r-1}dx-\int_0^{\mu} \left[\sum_{k=0}^n(-1)^k\binom{n}{k}x^k\right]x^rdx\nonumber\\
&=\sum_{k=0}^n \frac{1}{r+k}(-1)^k\tbinom{n}{k}\mu^{k+r+1}-\sum_{k=0}^n \frac{1}{r+k+1}(-1)^k\tbinom{n}{k}\mu^{k+r+1}\nonumber\\
&=\left[\sum_{k=0}^n \frac{1}{r+k}(-1)^k\tbinom{n}{k}\mu^k-\sum_{k=0}^n \frac{1}{r+k+1}(-1)^k\tbinom{n}{k}\mu^k\right]\mu^{r+1}.
\end{align}

By the binomial theorem, we have
\begin{align*}
\sum_{k=0}^n(-1)^k\binom{n}{k}\mu^kx^{k-n}=\left(\frac{1}{x}-\mu\right)^n.
\end{align*}
It follows that
\begin{align*}
\sum_{k=0}^n(-1)^k\binom{n}{k}\mu^kx^{r+k-1}&
=x^{n+r-1}\sum_{k=0}^n(-1)^k\binom{n}{k}\mu^kx^{k-n}\nonumber\\
&=x^{n+r-1}\left(\frac{1}{x}-\mu\right)^n,
\end{align*}
and thus
\begin{align*}
\sum_{k=0}^n \frac{1}{r+k}(-1)^k\tbinom{n}{k}\mu^k
&=\int_0^1\sum_{k=0}^n(-1)^k\binom{n}{k}\mu^kx^{r+k-1}dx\\
&=\int_0^1x^{n+r-1}\left(\frac{1}{x}-\mu\right)^ndx.
\end{align*}
Similarly, we have
 \begin{align*}
\sum_{k=0}^n \frac{1}{r+k+1}(-1)^k\tbinom{n}{k}\mu^k
&=\int_0^1\sum_{k=0}^n(-1)^k\binom{n}{k}\mu^kx^{r+k}dx\\
&=\int_0^1x^{n+r}\left(\frac{1}{x}-\mu\right)^ndx.
\end{align*}
Hence
\begin{align}\label{2.11}
&\sum_{k=0}^n \frac{1}{r+k}(-1)^k\tbinom{n}{k}\mu^k-\sum_{k=0}^n \frac{1}{r+k+1}(-1)^k\tbinom{n}{k}\mu^k\nonumber\\
&=\int_0^1\left(x^{n+r-1}-x^{n+r}\right)\left(\frac{1}{x}-\mu\right)^ndx\nonumber\\
&=\int_0^1x^{r-1}(1-x)\left(1-\mu x\right)^ndx\nonumber\\
&=\left.\frac{x^r}{r}(1-\mu x)^{n+1}\right|_0^1\nonumber\\
&=\frac{1}{r}(1-\mu)^{n+1},
\end{align}
where we used the fact that $\mu=\frac{r}{n+r+1}$ and
\begin{align*}
\frac{d\ \left(\frac{x^r}{r}(1-\mu x)^{n+1}\right)}{dx}
&=x^{r-1}(1-\mu x)^{n+1}-(n+1)\frac{x^r}{r}(1-\mu x)^n\cdot \mu\\
&=x^{r-1}(1-\mu x)^n\left[1-\mu x-\frac{n+1}{r}\cdot \mu x\right]\\
&=x^{r-1}(1-\mu x)^n\left[1-\mu x\left(1+\frac{n+1}{r}\right)\right]\\
&=x^{r-1}(1-x)\left(1-\mu x\right)^n.
\end{align*}

By (\ref{2.9})-(\ref{2.11}), we obtain (\ref{2.8}). The proof is complete.\hfill\fbox

\bigskip

\noindent {\bf Proof of Theorem \ref{thm}.}  By (\ref{1.1-b}), we have
\begin{align}\label{2.15}
P\left(NB(r,p)\leq \frac{r(1-p)}{p}\right)=\sum_{k=0}^{\left[\frac{r(1-p)}{p}\right]}\binom{k+r-1}{k}p^r(1-p)^k.
\end{align}

If $\left [\frac{r(1-p)}{p}\right]=n$ for some non-negative integer $n$,  then we have
$p \in \left(\frac{r}{n+r+1},\frac{r}{n+r}\right]$. Now, (\ref{2.15}) becomes
\begin{align}\label{2.16}
P\left(NB(r,p)\leq \frac{r(1-p)}{p}\right)=\sum_{k=0}^n\binom{k+r-1}{k}p^r(1-p)^k.
\end{align}

For any $k=0,\ldots,n$, define a function   $f_k(p):=p^r(1-p)^k, p \in \left(\frac{r}{n+r+1},\frac{r}{n+r}\right]$. For any $p \in \left(\frac{r}{n+r+1},\frac{r}{n+r}\right)$, we have
\begin{align*}
\frac{d\ln f_k(p)}{dp}&=\frac{r}{p}-\frac{k}{1-p}=\frac{1}{p(1-p)}\left[r-(r+k)p\right]\\
&>\frac{1}{p(1-p)}\left[r-(r+k)\frac{r}{n+r}\right]\\
&=\frac{r}{p(1-p)}\left[1-\frac{k+r}{n+r}\right]\\
&\geq 0,
\end{align*}
which implies that $\ln f_k(p)$ is strictly increasing on $\left(\frac{r}{n+r+1},\frac{r}{n+r}\right]$ and thus $f_k(p)$ is strictly increasing on $\left(\frac{r}{n+r+1},\frac{r}{n+r}\right]$. Then by (\ref{2.16}), we get that
\begin{align}\label{2.17}
\inf_{p \in \left(\frac{r}{n+r+1},\frac{r}{n+r}\right]}P\left(NB(r,p)\leq \frac{r(1-p)}{p}\right)=\left(\frac{r}{n+r+1}\right)^r
\sum_{k=0}^n\binom{k+r-1}{k}\left(1-\frac{r}{n+r+1}\right)^k.\end{align}

Define a sequence
$$
a_r(n):=\left(\frac{r}{n+r+1}\right)^r
\sum_{k=0}^n\binom{k+r-1}{k}\left(1-\frac{r}{n+r+1}\right)^k,\ n=0,1,\ldots.
$$
Obviously, we have
$$
a_r(0)=\left(\frac{r}{r+1}\right)^r.
$$
Then in order to finish the proof, it is enough to show that
$$
a_r(n)\geq a_r(0),\ \forall n=1,2,\ldots.
$$
In fact, we will prove the following stronger result:
\begin{align}\label{2.18}
a_r(n+1)-a_r(n)>0,\ \forall n=0,1,\ldots.
\end{align}

By Lemma 2.1, we have
\begin{align}\label{2.19}
a_r(n)=\frac{1}{\mathbf{B}(n+1,r)} \int_0^{\frac{r}{n+1+r}} x^{r-1}(1-x)^ndx.
\end{align}
By (\ref{2.19}), (\ref{2.0}) and (\ref{2.0-b}), we have
\begin{align}\label{2.20}
& (\ref{2.18})\ \mbox{holds}\nonumber\\
&\Leftrightarrow\frac{\Gamma(n+1+r)}{\Gamma(r)\Gamma(n+1)} \int_0^{\frac{r}{n+1+r}} x^{r-1}(1-x)^ndx < \frac{\Gamma(n+2+r)}{\Gamma(r)\Gamma(n+2)}\int_0^{\frac{r}{n+2+r}} x^{r-1}(1-x)^{n+1}dx\nonumber\\
&\Leftrightarrow\frac{n+1}{n+1+r} \int_0^{\frac{r}{n+1+r}} x^{r-1}(1-x)^ndx <\int_0^{\frac{r}{n+2+r}} x^{r-1}(1-x)^{n+1}dx\nonumber\\
&\Leftrightarrow\int_0^{\frac{r}{n+1+r}} x^{r-1}(1-x)^n\mathrm{d}x - \frac{r}{n+1+r} \int_0^{\frac{r}{n+1+r}} x^{r-1}(1-x)^ndx <\int_0^{\frac{r}{n+2+r}} x^{r-1}(1-x)^{n+1}dx.
\end{align}

By the change of variable $u=\frac{1}{x}$, we get
\begin{align*}
&\int_0^{\frac{r}{n+1+r}} x^{r-1}(1-x)^ndx=\int_{\frac{n+1+r}{r}}^{+\infty} (u-1)^{n} u^{-n-r-1}du,\\
&\int_0^{\frac{r}{n+2+r}} x^{r-1}(1-x)^{n+1}dx=\int_{\frac{n+2+r}{r}}^{+\infty} (u-1)^{n+1} u^{-n-r-2}du.
\end{align*}
Then, by (\ref{2.20}), we have
\begin{align}\label{2.21}
& (\ref{2.18})\ \mbox{holds}\nonumber\\
&\Leftrightarrow \int_{\frac{n+1+r}{r}}^{+\infty} (u-1)^{n} u^{-n-r-1}du- \frac{r}{n+1+r} \int_{\frac{n+1+r}{r}}^{+\infty} (u-1)^{n} u^{-n-r-1}du
<\int_{\frac{n+2+r}{r}}^{+\infty} (u-1)^{n+1} u^{-n-r-2}du\nonumber\\
&\Leftrightarrow\int_{\frac{n+1+r}{r}}^{+\infty} (u-1)^{n} u^{-n-r-1}du-
\left(\int_{\frac{n+1+r}{r}}^{+\infty}-\int_{\frac{n+1+r}{r}}^{\frac{n+2+r}{r}}\right)
(u-1)^{n+1} u^{-n-r-2}du\nonumber\\
&\quad\quad <\frac{r}{n+1+r} \int_{\frac{n+1+r}{r}}^{+\infty} (u-1)^{n} u^{-n-r-1}du\nonumber\\
&\Leftrightarrow\int_{\frac{n+1+r}{r}}^{\frac{n+2+r}{r}} (u-1)^{n+1}u^{-n-r-2}du< \int_{\frac{n+1+r}{r}}^{+\infty} (u-1)^n u^{-n-r-2}\left(\frac{ru}{n+1+r}-1\right)du.
\end{align}

Define a function
$$
g(u):=(u-1)^{n+1}u^{-n-r-2},\  u\in \left[\frac{n+1+r}{r},\frac{n+2+r}{r}\right].
$$
For any $u\in \left[\frac{n+1+r}{r},\frac{n+2+r}{r}\right]$, we have
\begin{align*}
g'(u)&=(u-1)^n u^{-n-r-3}[n+r+2-(r+1)u]\\
&<(u-1)^n u^{-n-r-3}\left[n+r+2-(r+1)\frac{n+1+r}{r}\right]\\
&=-\frac{(n+1)(u-1)^n u^{-n-r-3}}{r}\\
&<0,
\end{align*}
which implies that $g$ is  a strictly decreasing function and thus
\begin{align}\label{2.22}
\int_{\frac{n+1+r}{r}}^{\frac{n+2+r}{r}} (u-1)^{n+1}u^{-n-r-2}du
& <\frac{1}{r} \left(\frac{n+1+r}{r}-1\right)^{n+1} \left(\frac{n+1+r}{r}\right)^{-n-r-2}\nonumber\\
&=\frac{r^r(n+1)^{n+1}}{(n+r+1)^{n+r+2}}.
\end{align}

By  (\ref{2.22}) and Lemma 2.2, we obtain that the inequality in (\ref{2.21}) holds and thus (\ref{2.18}) holds.  The proof is complete. \hfill\fbox

\bigskip

\noindent {\bf\large Acknowledgments}\quad  This work was supported by the National Natural Science Foundation of China (12171335) and the Science Development Project of Sichuan University (2020SCUNL201).

\end{document}